\documentclass[12pt]{amsart} 
\usepackage{amsmath, amsthm, amscd, amsfonts,amssymb,graphicx,enumerate}
\numberwithin{equation}{section}
\makeatletter
\@namedef{subjclassname@2010}{%
  \textup{2010} Mathematics Subject Classification}
\makeatother

\newcommand{\teq}{\arabic{section}.\arabic{equation}}
\newcommand{\teql}{\Alph{section}.\arabic{equation}}

\newcommand{\sqr}[2]{{\vcenter{\vbox{\hrule height.#2pt\hbox{\vrule width.#2pt
height#1pt \kern#1pt\vrule width.#2pt}\hrule height.#2pt}}}}

\newcounter{eqcount}

\newenvironment{edesc}{\refstepcounter{equation}\begin{enumerate}}%
{\end{enumerate}}
\newenvironment{triv}{\refstepcounter{equation}\begin{list}%
{{\hbox{\rm(\teq)\ }}} \item }{\end{list}}

\newcommand{\ring}[1]{{\mathbb #1}}
\newcommand\bZ{{\ring{Z}}}

\newcommand\bC{{\ring{C}}} 
 \newcommand\bQ{{\ring{Q}}}

\newcommand{\csp}[1]{{\mathbb #1}}
\newcommand{\tsp}[1]{{\mathcal #1}}

\newcommand{\prP}{\csp{P}}
\newcommand{\afA}{\csp{A}}

\newcommand{\sO}{{\tsp{O}}}

\renewcommand{\ni}{\texto{Ni}}

\newcommand{\eql}[2]{{\rm (\ref{#1}\ref{#2})}} 

\newcommand{\vect}[1]{{\pmb #1}}

\newcommand{\bp}{{\vect{p}}} 
 
 \newcommand{\bz}{{\vect{z}}}

\newcommand{\row}[2]{{#1_1,\ldots,#1_{#2}}}

\newcommand{\smatrix}[4]{{\big(\begin{array}{cc}
\!\lower2pt\hbox{$\scriptstyle#1$} &\lower2pt\hbox{$\scriptstyle#2$}\!
\\\! \raise2pt\hbox{$\scriptstyle#3$} &\raise2pt\hbox{$\scriptstyle#4$}
\!\end{array}\big)}}

\newcommand{\texto}[1]{{\textr{#1}}}
 
 \newcommand{\ind}{\texto{ind}}
 \newcommand{\PGL}{\texto{PGL}}
 
 \renewcommand{\ni}{\texto{Ni}}
 
\newcommand{\textr}[1]{{\text{\rm #1}}}
\newcommand{\tr}{\textr{tr}} 
\newcommand{\abs}{\textr{abs}}  
 
 \newcommand{\inn}{\textr{in}}

\newcommand{\textb}[1]{{\text{\bf #1}}}
\newcommand{\bfC}{{\textb{C}}}
\newcommand{\longmapright}[2]{\smash{\mathop{\hbox to
#2pt{\rightarrowfill}}\limits^{#1}}}
\newcommand{\longmapleft}[2]{\smash{\mathop{\hbox to
#2pt{\leftarrowfill}}\limits^{#1}}}

\newcommand{\mapright}[1]{\smash{\mathop{\longrightarrow}\limits^{#1}}}

\newcommand{\np}{{+}}   \newcommand{\nm}{{-}}

\newcommand{\lrang}[1]{{\langle #1\rangle}}

\newcommand{\eqdef}{\stackrel{\text{\rm def}}{=}}

\newfont{\sevenrm}{cmr7}
\newfont{\bsevenrm}{cmbx7}
\newfont{\mathseven}{cmsy7}
\newfont{\bigmath}{cmsy10 scaled 1200}
\newfont{\fiverm}{cmr5}
\newfont{\bfiverm}{cmbx5}
\newfont{\hel}{cmbx10 scaled 1200}
\newfont{\eu}{eufb10}
\newfont{\sseu}{eufm5}
\newfont{\seu}{eufm7}
\newfont{\Cal}{cmmib10}
\newfont{\sCal}{cmmib7}
\newfont{\zch}{eusb10}

\theoremstyle{plain}
\newtheorem{thm}{Theorem}[section] 
\newtheorem{lem}[thm]{Lemma}

\newtheorem{prop}[thm]{Proposition}

\theoremstyle{definition}

\newtheorem{exmp}[thm]{Example}

\newtheorem{prob}[thm]{Problem}

\theoremstyle{remark}
\newtheorem{rem}[thm]{Remark}

\newcommand{\xs}{\times^s\!}

\def\pic #1 by #2 (#3){\vbox to #2{\hrule width 
#1 height 0pt depth 0pt\vfill\special{picture #3}}}
\def\scaledpicture#1
by #2 (#3 scaled #4){{\dimen0=#1 \dimen1=#2\divide\dimen0 by
 1000 \multiply\dimen0 by #4\divide\dimen1 by 1000 \multiply\dimen1 
by #4\pic \dimen0 by \dimen1 (#3 scaled #4)}}

\newcommand{\comm}[1]{{}}
\renewcommand{\phi}{\varphi}

\newcommand{\C}{{\text{\rm C}}}

\newcommand{\psigma}{{\pmb \sigma}}

\newcommand{\geng}{{{\text{\bf g}}}}

\begin{document}
\baselineskip=17pt

\title[Schinzel covers]{Schinzel's Problem: \\ \small{Imprimitive covers and the monodromy method}}

\author[M.~D.~Fried]{Michael
D.~Fried}
\address{Emeritus, UC Irvine \\ 3547 Prestwick Rd, Billings MT 59101}
\email{mfried@math.uci.edu}
\author[I.~Gusi\'c]{Ivica~Gusi\'c}
\address{faculty -- FKIT, Univ. of Zagreb  \\ Marulicev trg 19, Zagreb, Croatia}
\email{igusic@fkit.hr}

\date{}

\begin{abstract} There are now many successful uses of the {\sl monodromy method\/} for applying R(iemann)'s E(xistence) T(heorem) to describe solutions to problems on algebraic equations. Schinzel's original problem was to describe expressions $f(x)-g(y)$, with $f,g\in \bC[x]$ nonconstant, that are reducible. We call $(f,g)$ a {\sl Schinzel pair \/} if this happens {\sl  nontrivially\/} (see \eqref{trivfactors}).  

When $f$ is indecomposable \cite{Fr73} solved Schinzel's problem as a corollary. \cite{Fr10} revisits this to expand on many papers affected by the method: especially to circumvent using only covers with primitive monodromy group. We here take the next step  to consider the problem left by R.~Avanzi and U.~Zannier \cite{AZ03}, and the 2nd author \cite{Gu10}. Consider those  $f$ for which there is a $g=\alpha\circ f$, with $\alpha\in \PGL_2(\bC)$, satisfying an essential condition for possible Schinzel pairs: The Galois closure of the covers $f,g: \prP^1_x\to \prP^1_z$ are the same. Then, from those find $(f,g)$ that are Schinzel pairs.  
\end{abstract}




\thanks{}
\maketitle

\tableofcontents

\section{Schinzel's Problem and our particular case} \

For $f,g\in \bC[x]$, Schinzel's problem was to describe those cases when 
\begin{triv} \label{SchinzProb} \label{factors} $f(x)-g(y)$ factors nontrivially as a polynomial in two variables.  \end{triv} The topic is in \cite{Sc71};  \cite{Fr10} has many relevant references. With $K$ a number field, let $\sO_K$ be its ring of integers, $\bp$  a prime ideal of $\sO_K$, and $\sO_K/\bp$ its residue class field. Davenport's problem considered when, nontrivially,  
\begin{triv} \label{DavProb} the ranges of $f$ and $g$ are identical on almost all  $\sO_K/\bp$. \end{triv} 

The most trivial cases are where $g(x)=f(ax+b)$ for some $a,b\in \bar \bQ$, the algebraic numbers. When, $K=\bQ$, mostly that  relation forces $a,b\in K$. For example. this holds when $f$ is {\sl indecomposable\/} (not a composite of lower degree polynomials). With the indecomposable assumption, solutions to Davenport's and Schinzel's problems were essentially the same and solved (\cite[Thm.~1]{Fr73} and \cite[thm.~4.1]{Fr10}). 

Cases where $a,b$ aren't in $K$ are important to Davenport's problem, but not to Schinzel's. Though Schinzel's problem is our main concentration, in \S\ref{equivReps} the indecomposable case reappears in   Prob.~\ref{GusicProb},  our case of Schinzel's problem. The dihedral group, $D_n$, with $n$ even, the example of \S\ref{dihEx}, will aid a reader unaccustomed to branch cycles. Compare our goals with the \S\ref{GusicConj} conjecture.

\subsection{Branch cycles}  \label{branchCycles} We start by assuming $f=f_1\circ f_2$, and $\deg(f_i)>1$, $i=1,2$:  $f$ {\sl decomposes}. For Schinzel's Problem \eqref{SchinzProb} consider these extensions of what is a trivial relation between $f$ and $g$ (allowing a switch of $f$ and $g$). 
\begin{edesc} \label{trivfactors} \item  \label{trivfactorsa} Composition reducibility:  $ f_1(x) - g(y)$  factors.  
\item \label{trivfactorsb} A particular case of composition reducibility:  $g=f_1\circ g_2$, 
\end{edesc} \cite[Def.~2.1]{Fr87} calls an example of \eqref{factors} {\sl newly reducible\/} -- nontriviality for Schinzel's Problem --  if composite reducibility \eql{trivfactors}{trivfactorsa} does not hold. We call the corresponding $(f,g)$ a {\sl Schinzel pair\/}. 

We here consider the problem left by R.~Avanzi and U.~Zannier \cite{AZ03}, and the 2nd author \cite{Gu10}. Consider those  $f$ for which there is a $g=\alpha\circ f$, with $\alpha\in \PGL_2(\bC)$, satisfying an essential condition for possible Schinzel pairs: The Galois closure of the covers $f,g: \prP^1_x\to \prP^1_z$ are the same. Then, from those find $(f,g)$ that are Schinzel pairs.

Let  $\prP^1_z$ be the Riemann sphere, uniformized by the variable $z$.  Any rational function $f\in \bC(x)$ gives an analytic map -- a {\sl cover\/} -- $\prP^1_x\to \prP^1_z$. 
If the degree of $f$ is $n$, then {\sl branch points\/} of $f$ are the values of $z$ over which there are fewer than $n$ distinct points. For example, $z=\infty$ is a branch point of any polynomial $f\in \bC[x]$ with $\deg(f)>1$, because only $\infty$ lies over $\infty$. We denote the branch points of $f$ by $\bz_f=\{\row z r\}$. 

Refer to  $f_X: X\to \prP^1_z$, a compact Riemann surface cover, as Galois if the automorphisms that commute with $f_X$ have cardinality   $\deg(f_X)$. We often simplify $f_X$ to $f$ if there will be no misunderstanding. The Galois closure of $f$ is the smallest Galois cover, $\hat f: \hat X \to \prP^1_z$, that factors through $f$. It always exists. The group of automorphisms, $G_f$, of $\hat X$ commuting with $\hat f$ is the (geometric) monodromy group of $f$. 

The Galois correspondence associates to the cover $f_X$  a (faithful) coset (or permutation) representation  $T_f: G_f \to S_n$. We label a subgroup (up to conjugation by $G_f$) defining the cosets as $G(T_f,1)$. These are the elements of $G_f$ that fix the integer 1 in the representation $T_f$.  Similarly, any cover $f': X' \to \prP^1_z$ through which $\hat f$ factors corresponds to a coset representation (possibly not faithful) of $G_f$.

Whatever the branch points $\bz$, for any cover  $f : X \to \prP^1_z$ of compact Riemann surfaces, these produce conjugacy classes $\bfC={\row \C r}$ in the geometric monodromy $G_f\le S_n$. Denote  $\prP^1_z\setminus \{\bz\}$ by $ U_\bz$. \cite[\S5.3.2]{Fr10} explains using  {\sl classical generators\/} of the fundamental group of  $U_\bz$. These figure in  why you can select respective representatives $\sigma_i\in \C_i$, $i=1,\dots,r$,   to have these properties: 

\begin{edesc} \label{cycleconds} 
\item \label{cycleconds1} {\sl Generation}:  $\lrang{\sigma_i | i=1,\dots, r} = G_f\eqdef G \le  S_n$; and 
\item \label{cycleconds2} {\sl  Product-one}: $\sigma_1 \cdots  \sigma_r=1$.  \end{edesc} For fixed $\bfC$, the set of $\psigma$ satisfying \eqref{cycleconds} is the  {\sl Nielsen class}, $\ni(G,\bfC)$, of $(G,\bfC)$. Equivalences on Nielsen classes correspond to equivalences between covers. 

To get started we need only one: {\sl Absolute equivalence}. That means you mod out on Nielsen classes by the action of the subgroup of $S_n$, $N_{S_n}(G,\bfC)$, that normalizes $G$ and permutes (with multiplicity) the conjugacy classes in $\bfC$. The absolute equivalence class of $\psigma\in \ni(G,\bfC)$ 
is $$\{\alpha\psigma\alpha^{-1}\mid \alpha \in N_{S_n}(G,\bfC)\}.$$  Denote these  equivalence classes, running over $\psigma\in  \ni(G,\bfC)$ by $\ni(G,\bfC)^{\abs}$. 

The index, $\ind(\sigma)$, of a permutation $\sigma \in S_n$ is just $n$ minus the number of disjoint cycles in the permutation. Example: an $n$-cycle in $S_n$ has index $n\nm 1$, and an involution has index equal to the number of disjoint 2-cycles in it. The {\sl genus\/}, $\geng_X$ of $X$ given by $f_X$, with branch cycles in a given Nielsen class is well defined. The Riemann-Hurwitz formula says:  

\begin{equation} \label{RH}     2(n + \geng_X- 1) = \sum_{i=1}^r \ind(\sigma_i). \end{equation}

Two covers $f_i: X_i\to \prP^1_z$ are in the same {\sl absolute class\/} if there is a continuous (1-1) map $\psi: X_1\to X_2$ so that $f_1=f_2\circ \psi$. 

Further, the disjoint cycles of $\sigma_i$ correspond to  points of $X$ lying over $z_i$. A disjoint cycle length  is the ramification index of the point over $z_i$. An $r$-tuple, $\psigma$, satisfying \eqref{cycleconds} is a {\sl branch cycle description\/} of $f$. 
\cite[App.~A]{Fr10} explains {\sl classical generators\/} of the fundamental group of $U_\bz$ and how from them you get the following. 

\begin{prop} \label{theCorr} There is a 1-1 correspondence between elements of $\ni(G,\bfC)^{\abs}$ and absolute equivalence classes of covers  $f:X\to \prP^1_z$ in the Nielsen class, with any fixed set of $r$ distinct branch points $\bz$. 
\end{prop} 

We refer to Prop.~\ref{theCorr} as R(iemann's)E(xistence)T(heorem) or RET. \S\ref{procseta} uses  special  classical generators that work for our particular problem. 

\subsection{Reduced Galois-equivalence} 
Denote the functions $x\mapsto ax+b$, $a\in \bC^*$, $b\in \bC$, by $\afA(\bC)$. 
If $f,g\in \bC[x]$, and $g(x)=\alpha\circ f\circ \beta(x)$ (resp.~$f\circ \beta$), $\alpha, \beta \in \afA(\bC)$, we say  $f$ and $g$ are {\sl reduced  (resp.~affine) equivalent\/}. Call $f\in \bC[x]$ {\sl cyclic\/} if $f$ is reduced equivalent to $x^{\deg(f)}$. Consider the following for $(f,g)$ reduced, but not affine equivalent.   

\begin{edesc}  \label{redClos}  \item  \label{redClosa}  $f$ is not  cyclic  and  $f,g:\prP^1_x\to \prP^1_z$ have the same Galois closures.  
\item  \label{redClosb}  $f$ is not a composition of some polynomial with a non-trivial cyclic polynomial and $f(x)-g(y)$ is newly reducible. 
\end{edesc} 

We say the polynomials $f$ and $g$ satisfying \eql{redClos}{redClosa} are reduced {\sl Galois-equivalent}. With slight modification, the name makes sense for any pair of covers $f:X_f\to \prP^1_z$, $g: X_f\to \prP^1_z$, If they have the same Galois closure covers.  As in \S\ref{branchCycles} let $G_f(T_f,1)$ and $G_f(T_g,1)$ be the subgroups of $G_f$ corresponding to the covers $f$ and $g$. 

For  a cover represented by a non-cyclic polynomial, there is a unique branch cycle, $\sigma_\infty$ (attached to $z=\infty$), that has exactly one disjoint cycle (of length $n$).  
 
\begin{prop} \label{Gusic} If either of \eqref{redClos} hold, then translating $f$ by a constant, we may assume  $a=\zeta_v=e^{2\pi i/v}$, $v\ne 1$, and that $g=\zeta_v f$. Then, $a$ acts as a permutation $u_a$  of the finite branch points of $f$.   

If \eql{redClos}{redClosa} holds, then $z\mapsto az+b$ gives a cyclic cover $\mu: \prP^1_z\to \prP^1_{u}$ with group $\lrang{a^*}=\bZ/v$ where the following holds.  The composite covers $\mu\circ \hat f$ and $\mu\circ \hat g$ are also the same and Galois. If $\sigma_\infty^*\in G_{\mu \circ \hat f}$ is a branch cycle over $\infty$ for $\mu \circ \hat f$, then we can take its natural image in $\lrang{a^*}$  to be $a^*$, and $\sigma_\infty =(\sigma^*_\infty)^{v}$. 

Denote conjugation by $\sigma_\infty^*$ by $c_{\text{\rm AZ}}$. It has trivial action on $\sigma_\infty$ and no element of $S_n$ represents $c_{\text{\rm AZ}}$.   Up to conjugacy in $G_f$ we can choose  $c_{\text{\rm AZ}}$ to take $G_f(T_f,1)$ to $G_f(T_g,1)$. Identify $G_{\mu \circ \hat f}$ with the union of  $G_f$ cosets $\cup_{j=0}^{v\nm1} (\sigma_\infty^*)^jG_f$ (Rem.~\ref{CAZ}). 
\end{prop} 

\begin{proof}[About the proof of Prop.~\ref{Gusic}]  This is a special case of \cite[Prop.~7.28]{Fr10}. It stems from  \cite[Prop.~2]{Fr73}, which says -- under the newly reducible assumption -- that the Galois closures of $f$ and $g$ are the same. This general result has no dependence on the form of $f$ and $g$, except that their fiber product is newly reducible. Since their Galois closures are the same, their branch points are also identical. 

As $\sigma_\infty$ is a power of $\sigma_\infty^*$, $c_{\text{\rm AZ}}$ acts trivially on it. Since $\sigma_\infty^*$ normalizes $G_f$, it might be in $N_{S_n}(G_f)$. Yet, as it centralizes $\sigma_\infty$, it would have to be a power of $\sigma_\infty$ (the calculation \cite[Step 1, Proof Lem.~9]{Fr70});  contrary to it having order $v\cdot n$. 

The covers $f$ and $g$ correspond to  representations of $G$  on cosets of   $G_f(T_f,1)$ and $G_f(T_g,1)$. They are conjugate in $G_f$ if and only if $f$ and $g$ are absolutely equivalent covers:  the same as $f$ and $g$ being affine equivalent. By assumption they aren't.  So no element of $S_n$ represents $c_{\text{\rm AZ}}$. Choose the conjugates  $G_f(T_f,1)$ and $G_f(T_g,1)$ so that $\sigma_\infty^*$ conjugates one to the other. 
\end{proof} 

\begin{prob} \label{GusicProb} Characterize  branch cycles $\psigma$ (covers $f_X$) satisfying either of \eqref{redClos}.  For polynomials this includes $\geng_X=0$, but it makes sense without restricting  $\geng_X$. 
\end{prob} 

\begin{rem}[$c_{\text{\rm AZ}}$ leaves $\bfC$ invariant] \label{CAZ} The covers $f$ and $g$ in Prop.~\ref{Gusic} have the same Galois closures. So, it must be that $c_{\text{\rm AZ}}$ permutes the conjugacy classes in $\bfC$ -- preserving multiplicity -- just like the elements of $N_{S_n}(G,\bfC)$. 

Once we have identified the operator $c_{AZ}$ as in \S\ref{dihEx} or \S\ref{monMeth}, we can form $G_{\mu\circ \hat f}$ by taking a formal element $\sigma^*$, and forming the union of the left cosets of $G_f$. Multiplying coset elements comes from this formula for  $\sigma',\sigma''\in G_f$: 
$$(\sigma_\infty^*)^{j'}\sigma'(\sigma_\infty^*)^{j''}\sigma'' =
(\sigma_\infty^*)^{j'+j''}c_{AZ}^{-j''}(\sigma')\sigma''.$$ \end{rem} 
 
 \begin{exmp} \label{Chebychev} For $f\in \bC[x]$ and $a=-1$ in Prop.~\ref{Gusic}, $f$ and $-f$ define absolutely equivalent covers if and only if $f$ is affine equivalent to an odd $f^*$: $f^*(-x)=-f^*(x)$. That happens for odd degree in the general case of  \S\ref{dihEx}. Define the $n$th Chebychev polynomial, $T_n$, from  $T_n(\cos(\theta))$ being the real part of  $(e^{i\theta})^n=e^{ni\theta}$. For $n$ odd,  $T_n$ is odd from $(-e^{i\theta})^n=-(e^{i\theta})^n$.  \end{exmp} 
\subsection{Dihedral example, $D_n$,  $n$ even} \label{dihEx} Consider the semi-direct product, $$\bZ/n\xs A\eqdef \afA_n(A)\text{, with }A\le (\bZ/n)^*.$$ 
Regard it as the group of $2\times 2$ matrices: \begin{edesc} \label{afgroup} \item $ \bigl(\smatrix a b 01 \mid a\in A, b\in \bZ/n\bigr)$. With $A=\{\pm 1\}$, denote $\afA_n(A)$ by $D_n$. 
\item  \label{afgroupb}  Each element of $\afA_n(A)$ is a product $\smatrix 1 b 0 1 \smatrix a 0 0 1$.\end{edesc}  \cite[\S 7.2.1]{Fr10} (called \lq\lq Writing equations\rq\rq) gives the modern  -- but discusses the historical -- view of the  subgroups of  $\afA_n(A)$ playing the role of $G_{\mu \circ \hat f}$ in  Prop.~\ref{Gusic}. 
The set of involutions (order 2 elements) in $\afA_n(A)$ have the form $$I_n(A)=\{\smatrix a b 01 \mid a^2=1 (a\ne 1) \text{ and }b(a\np 1)=0.$$ 

\begin{lem} \label{auto} Assume $2|n\ge 4$. then, the distinct conjugacy classes,  $\C_{-1,0}$ and $\C_{-1,1}$, with reps.~$\sigma_1=\smatrix {-1} {1} 0 1$ and  $\sigma_2=\smatrix {-1} 0 0 1$,  comprise $I_n(\{\pm 1\})$. An   automorphism $c_n(\{\pm 1\})$, of $D_n$, is given by $$\smatrix 1 b 0 1 \mapsto \smatrix 1 b 0 1 \text{ and } \smatrix {-1} b 0 1 \mapsto \smatrix {-1} {b\nm1} 0 1,\  b\in \bZ/n. $$ 
\end{lem}  

The lemma follows easily by computation, with $c_n(\{\pm 1\})^2$ the same as  conjugation by $\sigma_\infty=\smatrix {-1} {-1} 0 1$. 
Now with $v=2$, $\zeta_v=-1$, we describe $D_n^*$ so it fits the conclusion of Prop.~\ref{Gusic} as $G_{\mu\circ f}$. Use $\sigma_i$, $i=1,2$, from Lem.~\ref{auto}. As generators $D_n^*$ has  $\sigma_1$,  and $\sigma_\infty^*$ satisfying these conditions: 
\begin{equation} \label{extProperty} \begin{array}{rl} (\sigma_\infty^*)^2 =&\sigma_\infty \text{ ($\sigma_\infty^*$ has order $2n$); and} \\
(\sigma_\infty^*)^{k}\smatrix{-1} 0 0 1(\sigma_\infty^*)^{-k} =&\smatrix{-1}  {-k}  0 1. \end{array}\end{equation}  
Denote the representation, from permutations in \eqref{chebyPer} by $T_f$. It comes from acting on (left) cosets of  $\lrang{\sigma_2}$. Another representation,  $T_g$, comes from cosets of $\lrang{\sigma_1}$.

The corresponding cover  -- given by a degree $n$ Chebychev polynomial --  appears in Prop.~\ref{Gusic} with $r=3$ and $A=\{\pm 1\}$. Assume $2| n\ge 4$. With the elements acting as permutations -- from the left -- on the integers $\{0,1,\dots, n\nm 1\} \mod n$: \begin{equation}\begin{array}{rl} \label{chebyPer} \sigma_1=& (1\,n)(2\,n\nm1)\cdots (\frac n 2\,\frac n 2\np1)\\
\sigma_2=&(1\,n\nm1)(2\,n\nm2)\cdots (\frac n 2\nm 1\,\frac n 2\np1) \\
\sigma_\infty\eqdef \sigma_3&=(1\,2\,\dots n\nm1\, n)^{-1} =\smatrix 1 {-1} 0 1. \end{array}\end{equation} 

We now apply RET (Prop.~\ref{theCorr}) to \eqref{chebyPer} to produce a polynomial pair $(f,g)$ with these properties for any even $n\ge 4$: 
\begin{edesc} \item all irreducible factors of $f(x)- g(y)$  have degree 2; but in this case
\item  \eqref{factors}  is newly reducible only for $n=4$.  \end{edesc} 
For finite branch points take any pair $(z',-z')$. For simplicity  we'll take $z'=1$. As in \S\ref{branchCycles}, $f$ (resp.~$g$) corresponds to the permutation representation $T_f$ (resp.~$T_g$). 

Respective indices of the $\sigma_i\,$s  in \eqref{chebyPer}, are $\frac n 2$, $\frac n {2} \nm 1$, and $n\nm 1$. Plug these into \eqref{RH} and conclude the genus of the cover -- call it $f$ -- is 0. Similarly, for a cover $g$ from $T_g$. Now use this characterization: $f: X \to \prP^1_z$ is absolutely equivalent to a polynomial cover if $X$ has genus 0, and precisely one point lies over $z=\infty$. 

RET and the Galois correspondence give the following. 
\begin{edesc} \label{chebBr} \item \label{chebBra} The irreducible factors of $f(x)-g(y)$ correspond one-one with the orbits of $G(T_g,1)$ in $T_f$ (on the cosets of $D_n(T_f,1)$), all length 2.
\item \label{chebBrb} The representation $T_{\mu\circ f}$ corresponding to $\mu\circ f: X \to \prP^1_u$, having monodromy $D_n^*$, is on the $2\cdot n$ cosets of $D_n^*(T_f,1)$. 
\item  \label{chebBrc} Composing $c_{\text{\rm AZ}}$  (as in \S\ref{equivReps}) in Prop.~\ref{Gusic} with $T_{\mu\circ f}$ is equivalent to $T_{\mu\circ g}$.  \end{edesc} 

Finally,  $D_n$ maps to $D_{\frac n 2}$, with a compatible representation on the cosets of $\lrang{\sigma_1}$. Then, $f$ is a composite of  degree 2 and $\frac n 2$ polynomials. When $\frac n 2=n'$ is odd, use Ex.~\ref{Chebychev} to see the replacement for  \eqref{chebBr}  (as in \cite[Lem.~7.4]{Fr10}) has a factor of degree 1, the rest of degree 2. So, unless $\deg(f)=4$, $f$ is not newly reducible. 

We use the principle \lq dragging a cover by its branch points\rq\ (\cite[\S 6.1]{Fr10}) to  producing a new cover from $f:X\to\prP^1_z$ with the same branch cycles, but finite branch points placed at any distinct points in $\bC$. We  require  $\zeta_v$ to permute the finite branch points. Example: for orbit condition \eqref{oneOrbitCond}, we may assume  $f$ has finite branch points $\zeta_v^j$, $j=1,\dots,v$. Then, $\mu\circ f$ has branch points 0, 1 and $\infty$. 
\begin{prob} \label{compBranch} As in \eql{chebBr}{chebBrc}, compute branch cycles, $\sigma_0^*, \sigma_1^*, \sigma_\infty^*$ for $\mu\circ f$. \end{prob}  

\begin{proof}[Hints for Prob.~\ref{compBranch}]   Take the branch cycle for $\infty$ as $\sigma_\infty^*$ using   \eqref{sigma*bcyc}. The shape of the branch cycle for  $\sigma_1^*$ is a product of $\frac{n\np n\nm 2}2=n\nm 1$   disjoint 2-cycles,  from juxtaposing contributions of $\sigma_1$ and $\sigma_2$; and $\sigma_0^*$ is a product of $n$ disjoint $v$-cycles. 

The Nielsen class $\ni(D_n,\bfC)^\abs$with the 3 conjugacy classes represented in \eqref{chebyPer} has  6 elements, indicated by the order of those conjugacy classes in a representing 3-tuple. This is common when $r=3$, but for $r\ge 4$, the braid group enters, as used in \cite[\S6.4]{Fr10}  to  dramatic effect. So, here inspection can produce the desired $\sigma_0^*,\sigma_1^*,\sigma_\infty^*$ satisfying generation and product-one in \eqref{cycleconds}. 

For, however,  a general polynomial $f$ only knowing that $u_a$ in Prop.~\ref{Gusic} permutes   
the branch points,  solving this problem requires {\sl classical generators\/} (as in \S\ref{monMeth}) for the covers $f$ and $\mu\circ f$,  and then a relation between these as in \S\ref{Ba02}. 
\end{proof} 

\subsection{The conjecture of \cite{Gu10}} \label{GusicConj} Our tentative conjecture is that \S\ref{dihEx} (with $n=4$) gives the only case of Schinzel pairs of the form $(f,\zeta_vf)$. As the argument of \cite[p.~47]{Fr70} shows, this  is true if and only if $\sigma_\infty$ generates a normal subgroup in $G$. 

Precisely: The conjugation $\sigma_i\sigma_\infty\sigma_i^{-1}=\sigma_\infty^k$  by a finite branch cycle implies $\sigma_i$ has  the same  index as multiplication by  $k\in (\bZ/n)^*$ on $\bZ/n$. Possibilities for a genus 0 cover (using \eqref{RH}) shows  $f$ is equivalent to a Chebychev (or cyclic) polynomial, with well understood branch cycles. Then, the Nielsen class -- according to \S\ref{dihEx} -- must be $\ni(D_4,\bfC)^\abs$ with $\bfC=\C_{-1,0}\cup \C_{-1,1}\cup\C_\infty$, as in Lem.~\ref{auto}. 

\section{The Group formulation of conditions \eqref{redClos}}  \label{monMeth} A conclusion from Prop.~\ref{Gusic} is that (as in the last Hint to Prob.~\ref{compBranch}): 
\begin{equation} \label{branchPtCond}  \bz_f=\zeta_v\bz_f=\bz_g.\end{equation}  
Consider  any polynomial $f$ assuming \eqref{branchPtCond}, and  one further condition: 
\begin{triv} \label{oneOrbitCond} $u_a$ has one orbit on finite branch points: $r-1=v$. \end{triv} 
\noindent That is, $\bz_f$ are the vertices of a  regular $v$-gon on a circle around the origin. 
\S\ref{procseta} sets up the procedure for computing branch cycles for  $\zeta_vf$ from those of $f$.  

\S\ref{procsetb} then characterizes possible branch cycles when you add \eql{redClos}{redClosa}, the Galois closure assumption for the pair $(f,\zeta_vf)$. \S\ref{manyCycle} notes that we can adjust the method to handle Prop.~\ref{Gusic} without condition \eqref{oneOrbitCond}.  

\subsection{The effect of $u_a$ on branch cycles when \eqref{branchPtCond} holds}  \label{procseta} 
Let $\afA_{r\nm 1}^0$ consist of all distinct $r\nm 1$-tuples in $\bC$. Assume \eqref{oneOrbitCond}. Then, given branch cycles for $f$ relative to classical generators of $\pi_1(U_{\bz_f}, z_0)$ it makes sense to compute branch cycles for $\zeta_vf$ relative to the same classical generators. 

Since we have assumed $z=0$ is not a branch point, we can use it as a basepoint, and the paths of App.~\ref{classPolyGens} -- where $r\nm1 =6$ -- listed as $\row {\bar \sigma} {r\nm1}, \bar \sigma_\infty$. 
We can compose any cover $f: X\to \prP^1_z$ with any element $\alpha\in \PGL_2(\bC)$. We make an increasing sequence of assumptions, starting with this: 
\begin{triv} \label{specAlpha} Suppose $\alpha(\bz_f)=\bz_f$: $\alpha$ permutes the branch points. \end{triv} We only do the next lemma for the case we use in the rest of the paper, and with the classical generators  $\row {\bar \sigma} r=\bar\psigma$ of App.~\ref{classPolyGens}. 

\begin{lem} We can explicitly compute the effect of $\alpha$ on an explicit set of classical generators to find branch cycles for $\alpha\circ f$ from branch cycles for $f$. Assume $\alpha$ is multiplication by $\zeta_v$ under assumption \eqref{oneOrbitCond} and $\psigma$ are branch cycles for $f$ relative to $\bar \psigma$ above. Then, relative to $\bar \psigma$, branch cycles for $\zeta_v f$ are 
\begin{equation} \label{oneCycle} (\sigma_2,\dots, \sigma_{r\nm 1},\sigma_1, \sigma_1^{-1}\sigma_r\sigma_1). \end{equation} \end{lem} 

\begin{proof}  Rotation through an angle of $2\pi/v$ sends $\bar\sigma_i$ to $\bar \sigma_i'=\bar \sigma_{i\np1}$, $i=1,\dots, r\nm 2$, and $\bar \sigma_{r\nm 1}$ to $\bar \sigma_{r\nm 1}'=\bar \sigma_1$. Similarly, $\bar \sigma_r$ (on the meridian halfway between $\bar \sigma_r$
and $\bar \sigma_1$) rotates to the meridian halfway between $\bar \sigma_1$ and $\bar \sigma_2$.  

Here's the deal! The branch cycles for $\alpha\circ f$ relative to  $\bar \psigma'$ are $\psigma$, the same as those for $f$ computed relative to the $\row {\bar \sigma} r$. Write $\row {\bar \sigma} r$ -- up to isotopy -- as words in $\row {\bar \sigma'} r$. Then,  plug $\psigma$ in to get the branch cycles for $\alpha\circ f$. To do that we only need express $\sigma_r$ by the following formula. Up to isotopy 
\begin{equation} \label{theRelation} \bar\sigma_1\bar \sigma_r'=\bar \sigma_r\bar \sigma_1. \end{equation} Explanation: The left side deforms on $U_{\bz_f}$ -- without moving $z_0$, or touching any points of the paths outside of $\bar \sigma_r,\bar \sigma_1,\bar \sigma_r'$ -- to a \lq\lq circle\rq\rq\ based at $z_0$ around $z_1$ and $\infty$. This is homotopic to a deformation of the right side of  \eqref{theRelation} that does the same. \end{proof} 

\subsection{Adding the Galois closure condition}  \label{procsetb} Suppose we have $(G,\bfC)$, with two  (faithful) permutation representations $T_i: G \to S_{n_i}$, $i=1,2$. (Our example will have $n_1=n_2=n$.) Then, we have two absolute Nielsen classes: $\ni(G,\bfC)^{\abs,i}$, $i=1,2$. Assume, too, we have ${}_i\psigma \in \ni(G,\bfC)^{\abs,i}$, representative classes and,  as in Prop.~\ref{theCorr}, these define covers $f_i: X_i\to \prP^1_z$, $i=1,2$, with branch points  $\bz$,  relative to specific classical generators. 

We must add {\sl inner equivalence\/} to {\sl absolute equivalence\/} on Nielsen classes (\S\ref{branchCycles}), $\ni(G,\bfC)^\inn\eqdef \ni(G,\bfC)/G$, to formulate the criterion that the $f_i\,$s have the same Galois closure covers. That is,  mod out by just  $G$ acting inside $N_{S_n}(G,\bfC)$. 

\cite[\S B.2.1]{Fr10} uses  examples to show how absolute and inner classes relate -- starting from  the canonical maps $\psi_{\inn,\abs}: \ni(G,\bfC)^\inn\to \ni(G,\bfC)^{\abs}$ -- to the main ingredient of \cite[Main Thm.]{FrV91}. The following is a natural addendum. 

\begin{prop}  \label{charSchinzel} The covers $f_1$ and $f_2$ have the same Galois closures if there exists $\psigma\in \ni(G,\bfC)^\inn$ for which $\psi_{\inn,\abs,i}(\psigma)={}_i\psigma$, $i=1,2$. The following characterizes  there being  a polynomial  $f$ in $\ni(G,\bfC)^\abs$ with $g=\zeta_vf$ satisfying  \eql{redClos}{redClosa} (Galois closure condition), with branch points $\bz$ satisfying one-orbit condition \eqref{oneOrbitCond}. \begin{edesc} \item $G$ has an automorphism, conjugation by $\sigma_\infty^*$, as in Prop.~\ref{Gusic}, with  
\item $\psigma\in \ni(G,\bfC)^\abs$  of genus 0 (a la \eqref{RH}),  $r\nm 1=v$, and  \eqref{oneCycle} holds. \end{edesc} 

There is an analog for more general orbits of $\zeta_v$. See \S\ref{manyCycle}. 
\end{prop} 

We note two points about the \S\ref{dihEx} example. 1st: We checked separately that we got reduciblity (for all $n$).  Then,  that it gave newly reducible, so a Schinzel pair (as in \eql{redClos}{redClosb}) just in the case $n=4$. Still, both came directly from branch cycles.  It is easy to generalize those conditions to apply to Prop.~\ref{charSchinzel}. 

2nd: In \S\ref{dihEx} conjugation by $\sigma_\infty^*$ permutes two distinct conjugacy classes. \S\ref{equivReps} shows we must have something like that to get Schinzel pairs. 

\subsection{Characterizing the $f$ in Prop.~\ref{Gusic} in general}  \label{manyCycle} 
Although more intricate, we can generalize \eqref{oneCycle} for any number of orbits for multiplication by $\zeta_v$ on branch points. It is possible, that with more than one orbit, we might have the origin as a branch point. We hope to complete the one orbit case of this paper in a later paper. There we will treat the generalization of \eqref{oneCycle}. 

 \subsection{Equivalent representations} \label{equivReps} We continue the 2nd observation at the end of \S\ref{procsetb}. Assume in Prop.~\ref{Gusic}  that $\sigma_\infty^*=\sigma^*$ satisfies the following condition: 
  \begin{triv} \label{autassume}  Conjugation, $c_{\sigma^*}$, by $\sigma^*$ preserves all conjugacy classes. \end{triv} 
\noindent Indeed, we aim for generality for future use. Assume a finite group $G$ has an outer automorphism  $\gamma$ (in place of $c_{\sigma^*}$) {\sl preserving classes\/}  -- the conclusion of \eqref{autassume}. 

Then, Prop.~\ref{primsystems} shows   $f$ could not possibly give new Schinzel pairs.  Applied to the conditions of Prop.~\ref{charSchinzel} it does produce a variables separated factorization $f(x)-g(y)$, but it  is not newly reducible: \eql{trivfactors}{trivfactorsa} holds.   
It still may contribute to Davenport's problem \eqref{DavProb} where, if the range values are assumed with the same multiplicities, the representations $T_f$ and $T_g$ satisfy the conclusion of Lem.~\ref{autoconj}.  

Applying any automorphism, $\gamma$,  to any permutation representation $T: G\to S_n$ sends it to another representation: $$T_{\gamma}: \sigma \mapsto T\circ \gamma(\sigma), \sigma\in G.$$ Denote  the stabilizer of an integer in $T$ by $G(T,1)$ and the number of fixed integers of $T(\sigma)$ by $\tr(T(\sigma))$: it's {\sl trace\/}.  
 
\begin{lem} \label{autoconj} Consider a representation $T: G \to S_n$. Assume $\gamma$ preserves classes. Then $\tr(T(\sigma))=\tr(T_\gamma(\sigma))$ for all $\sigma \in G$. \end{lem} 

\begin{proof} We are comparing the cosets of $G(T,1)$ fixed by $\sigma$ (multiplying on the left) with the cosets  fixed by $\gamma(\sigma)$. Since conjugation by $\gamma$ preserves the conjugacy class of $\sigma$: $\gamma(\sigma)=\sigma'\sigma (\sigma')^{-1}$ for some $\sigma'\in G$. The fixed cosets of $\sigma'\sigma (\sigma')^{-1}$ are the same as the fixed cosets of $\sigma$ on the conjugates of those cosets by $\sigma'$. But, if $T(\sigma')(1)=k$, then $\sigma'$ conjugates those cosets to the cosets of $G(T,k)$. Now, $\sigma$ fixes exactly the same number of $G(T,1)$ cosets as it fixes of $G(T,k)$  cosets. We are done. \end{proof}

Suppose we start with a fixed faithful transitive permutation representation $T_f$, coming from a cover of nonsingular curves $f: X_f\to Y$ (over $\bC$). Apply the Galois correspondence. It gives  a one-one correspondence between (nonsingular) covers  $f': X'\to Y$ through which $f$ factors, up to absolute equivalence, and groups $G(T,1)\le G'\le G$. Each $G'$  corresponds  to a system of imprimitivity of the permutation representation. This generalizes the notion of composition factors of a polynomial (or rational function).  

\begin{prop} \label{primsystems} Assume $\gamma$ and $T$ as above, with $g: X_g\to Y$ corresponding to $T_\gamma$. Then, the (normalization of the) fiber product $X_f\times_YX_g$ is reducible. This applies to the permutation representation $T'$ attached to any $G'$ with $G(T,1) < G' < G$. In particular, if $T$ is not primitive, then $X_f\times_YX_g$ is not newly reducible. 

So, if $(f,g)$ is a polynomial pair from Prop.~\ref{Gusic}, with $G=G_f$, $\gamma=c_{\sigma_\infty^*}$, then $f(x)-g(y)$ is reducible (as in \eqref{SchinzProb}), but not newly reducible. \end{prop}

\begin{proof} \cite[\S2.3]{Fr10} discusses Galois Theory and fiber products. Including that we naturally form the Galois closure of a degree $n$ cover from a component of the fiber product of the cover with itself, taken $n$ times.  Thus, the two topics go together: use of normalization (which for curves means the results are nonsingular); and how this generalizes the case of two polynomials $(f,g)$ as in the last statement. 

This paper's case (over the complexes) is  easier than in  \cite{Fr10}, over any characteristic zero field. The point is to have Galois theory  turn statements relating two covers into statements comparing two permutation representations. For example, consider this  statement:  $X_f\times_YX_g$ is reducible, which  \cite[\S 2.1]{Fr10} shows generalizes saying  \eqref{SchinzProb}.  The translation is that 
\begin{triv} \label{twoOrbits} $G_f(T_g,1)$ has more than one orbit in the representation $T_f$.\end{triv}  

This exactly generalizes \eql{chebBr}{chebBra} in \S\ref{dihEx}, except we computed directly that all orbits there had length 2 (for  $n > 2$). Here, a short argument from group theory applies:  \cite[Lem.~3]{Fr73} and assiduously redone in 
\cite[Rem.~4.3]{Fr10}, titled \lq\lq Davenport without $f$  indecomposable.\rq\rq\ It says \eqref{twoOrbits} follows from the weaker condition 
\begin{triv}  $\tr(T_f(\sigma))> 0$ if and only if $\tr(T_g(\sigma))>0$ for all $\sigma \in G$.\end{triv} 

Similarly, consider how we figured that only for $n=4$ would the \S\ref{dihEx} example be newly reducible. In our general case we assumed $T$ is not primitive. So, there is a representation $T'$ on the cosets of a group properly between $G(T,1)$ and $G$. According to Lem.~\ref{autoconj}, this produces the two representations $T'$ and $T'_\gamma$ to which we can apply the reducibility result above. We only need, in the last sentence, where $(f,g)$ are polynomials, that  $f$ decomposes under the hypotheses of Prop.~\ref{Gusic}. This is in the paragraphs above   \cite[Rem.~7.7, at the end of \S7.2.3]{Fr10} (called \lq\lq Ritt I\rq\rq) where we revamped how \cite{AZ03} treated the indecomposable case.  
\end{proof} 

\section{Searching for $(G,\bfC)$ that give Schinzel pairs} \label{sfactoring}  These short comments  suggest tools for dealing with what remains unsolved here, or on related problems. \S\ref{Ba02} and \S\ref{compCyc} are additions to the wreath product comments of \cite[\S7.2.4]{Fr10}. \S\ref{compCyc} focuses on our main case: $G_\mu=\bZ/v$. 

\subsection{Comments on \cite{Ba02}} \label{Ba02} Suppose we have any sequence of covers $$X \mapright{f }\prP^1_x \mapright{\mu} \prP^1_u.$$  \cite{BiFr86} and \cite{Tr93}  provide results for more general problems where the target of $f$ is not necessarily genus 0.  Simplifying, however,  for our special case is the work of \cite[Chap.~V]{Ba02}, called \lq\lq Nielsen graphs.\rq\rq\  

\begin{triv} \label{brcycGoals} From branch cycles for $\mu\circ f$ (relative to its base's classical generators), we can compute branch cycles for $f$. \end{triv}   In the other direction, branch cycles for $f$ and $\mu$ give information on branch cycles for $\mu \circ f$. We naturally identify the monodromy group $G_{\mu\circ f}$ with a subgroup of the {\sl wreath product\/} of $G_f$ and $G_\mu$, $G_f\wr G_\mu$, a completely general statement. 

For general $\mu$ of degree $v$, $G_f\wr G_\mu$ is naturally the semi-direct product $(G_f)^v\xs G_\mu$. Suppose $G_\mu\le S_v$. Denote the $i$th copy of $G_f$ in  $(G_f)^v$  by $G_{f,i}$, $i=1,\dots,v$. Then, here is the action of  $\gamma\in G_\mu$: $$(\row \sigma v)\in (G_f)^v \mapsto (\sigma_{(1)\gamma}, \dots ,\sigma_{(v)\gamma}).$$ It permutes the coordinates 
of  $(G_f)^v$ according to the permutation effect of $\gamma$. 

Suppose  $f$ and $\mu$ both are polynomial covers. According to \cite[Lem.~15]{Fr70}, $G_{\mu\circ f}$ will be the full wreath product under the following conditions. 
\begin{triv} \label{disjointness} The image of the finite branch points of $f$ under $\mu$ are all distinct and also distinct from the (finite) branch points of $\mu$. \end{triv}

If  the conditions of \cite[Lem.~15]{Fr70} don't hold, then $G_{\mu\circ f}$ may be a {\sl proper\/} subgroup of $G_f\wr G_\mu$, but still satisfying these conditions: 
\begin{triv} \label{wreathconds} $G_{\mu\circ f}$ maps surjectively onto  $G_\mu$, and its intersection with $(G_f)^{v}$  projects surjectively onto each $G_{f,i}$, $i=1,\dots,v$. \end{triv} 

\subsection{The case $G_\mu=\bZ/v$} \label{compCyc} To simplify notation we  use a superscript $*$-notation for elements in the group $G_{\mu\circ f}\eqdef G^*$. Our basic assumptions wll be the following: 
\begin{edesc} \label{fconds} \item \label{fcondsa} $\mu(z)=z^v$ is a cyclic cover of degree $v>1$; 
\item \label{fcondsb} $f$ is in a genus 0 Nielsen class, totally ramified over $z=\infty$; and 
\item \label{fcondsc}  the finite branch points of $f$ fall into $s$ orbits of (exact) length $v$ under multiplication by $e^{2\pi i/v}$.\end{edesc} 
From \eql{fconds}{fcondsc}, $\bfC$ has $r\nm 1 =s\cdot v$ conjugacy classes in it corresponding to finite branch point. From \eql{fconds}{fcondsa}, the branch points in each $e^{2\pi i/v}$ orbit go to the same value of $u$ under $\mu$. Finally, from \eql{fconds}{fcondsb}, $\mu\circ f$ totally ramifies over $\infty$, corresponding to a branch cycle $\sigma_\infty^*$ that has order $n\cdot v=n^*$. 

A description of branch cycles for the cover $\mu\circ f$ includes a branch cycle at $\infty$, given by an $n\cdot v$-cycle $\sigma^*$. We now write notation for $\sigma^*$. 
Identify $v$ copies of $\{1,\dots,n\}$ as $\{1_i,\dots,n_i\}$, the integers on which $G_{f,i}$ acts, $i=1,\dots, v$. With no loss, up to renaming the letters -- using that $(\sigma_\infty^*)^v=\sigma_\infty$ -- we can take $\sigma_\infty^*$ as  
\begin{equation} \label{sigma*bcyc} (1_1\,1_2\,\dots\,1_d\,2_1\,\dots\,2_d\,\dots \,n\nm1_1\,\dots \,n\nm1_d\,n_1\,\dots\,n_d).\end{equation}  Then, $\sigma_\infty$ generates the  intersection of $\lrang{\sigma_\infty^*}$ with $(G_f)^v$. 

In our situation, as  in Prop~\ref{Gusic}, the actual $G_{\mu \circ \hat f}\eqdef G_{f^*}$ is  the smallest subgroup of the full wreath product, $G_f\wr \bZ/v=(G_f)^v \xs \bZ/v$, satisfying  wreath conditions \eqref{wreathconds}. 

\subsection{New non-polynomial Schinzel pairs} Prop.~\ref{primsystems} produces general  fiber products of covers that may not have genus 0. With its extra hypothesis, however, these aren't newly reducible. To expand our understanding of Schinzel pairs  we might drop the condition they come from polynomials or even that they come from genus 0 covers.  Yet, they produce new fiber products,  from a pair of covers $f: X\to \prP^1_z$ and $g=\zeta_v f$, for which $f(x)-g(y)$ is newly reducible. 

Ex.~\ref{modular} uses genus 0 covers, given by rational functions, rather than polynomials. Again $v=2$, but $\zeta_2=-1$ has two orbits on  four finite branch points. 

\begin{exmp} \label{modular} Here $r=4$. Use the \lq dragging a cover by its branch points\rq\ principle of \S\ref{dihEx} to place the branch points at -1, -2, +2, +1  to correspond to branch cycles $(\sigma_1,\sigma_2,\sigma_2,\sigma_1)$ as given in \eqref{chebyPer}. The Nielsen class here contains the two conjugacy classes labeled $\C_{-1,0},\C_{-1,1}$, both twice, but it does not include an $n$-cycle. The group is still $D_n$;  the Galois closure has genus 1 (not 0 as in \S\ref{dihEx}). 

Many -- as a function of $n$ -- covers in the Nielsen class, correspond to different branch cycles. Yet, only two give a $g=-f$ with the same Galois closure. To be precise we must give classical generators replacing those of App.~\ref{classPolyGens}. They are almost the same \lq lolly-pop\rq\ paths from the origin through -1, -2, +2, +1, except, you can't allow the  lolly-pop  that passes around  -2 to go through -1. Instead, take a little blip to the right around -1 before continuing onto the rest of the lolly-pop. Similarly for the lolly-pop through +2, a little blip to the left around +1. 
\end{exmp} 

Now we suggest how to get new groups, but with covers of genus $>0$. 
\begin{prob} Extend the automorphism $c_{n}$ of Lem.~\ref{auto} to other subgroups of $\afA(n)$ to produce new, newly reducible fiber products  analogous to  $n=4$ of \S\ref{dihEx}. \end{prob} 

\begin{appendix} 
\section{Regular polygon classical generators} \label{classPolyGens} 
The paths, $\delta_i\sigma^*_i\delta_i^{-1}$ (including the subscript $r=\infty$, going around $\infty$ in Fig.~1 satisfy all the conditions of {\sl classical generators\/} based at $z_0=0$. Our notation is compatible with that of \cite[App.~B.1]{Fr10}, except we here use very regular paths, with punctures (except at $\infty$) arranged on a regular 6-gon. 

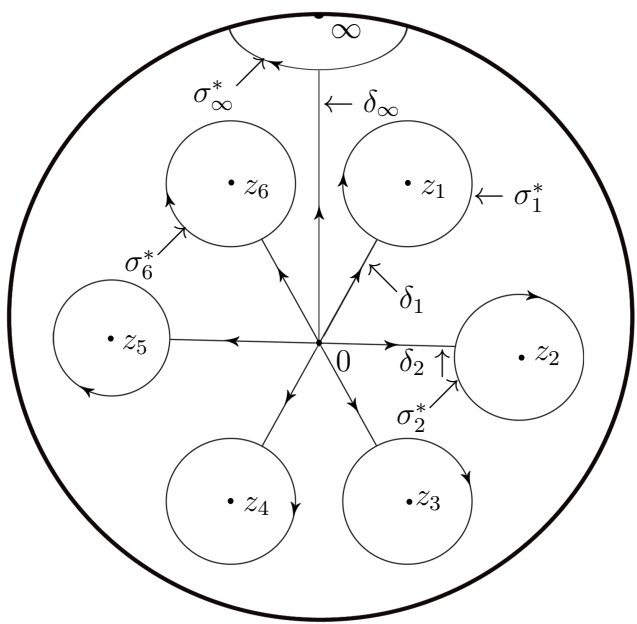
\begin{figure}[h]
\caption{$r=7$, with 6 branch points on a regular polygon }

\!\!\!\!\!\!\!\!$$
\hbox{\centerline{\includegraphics[scale=.6]{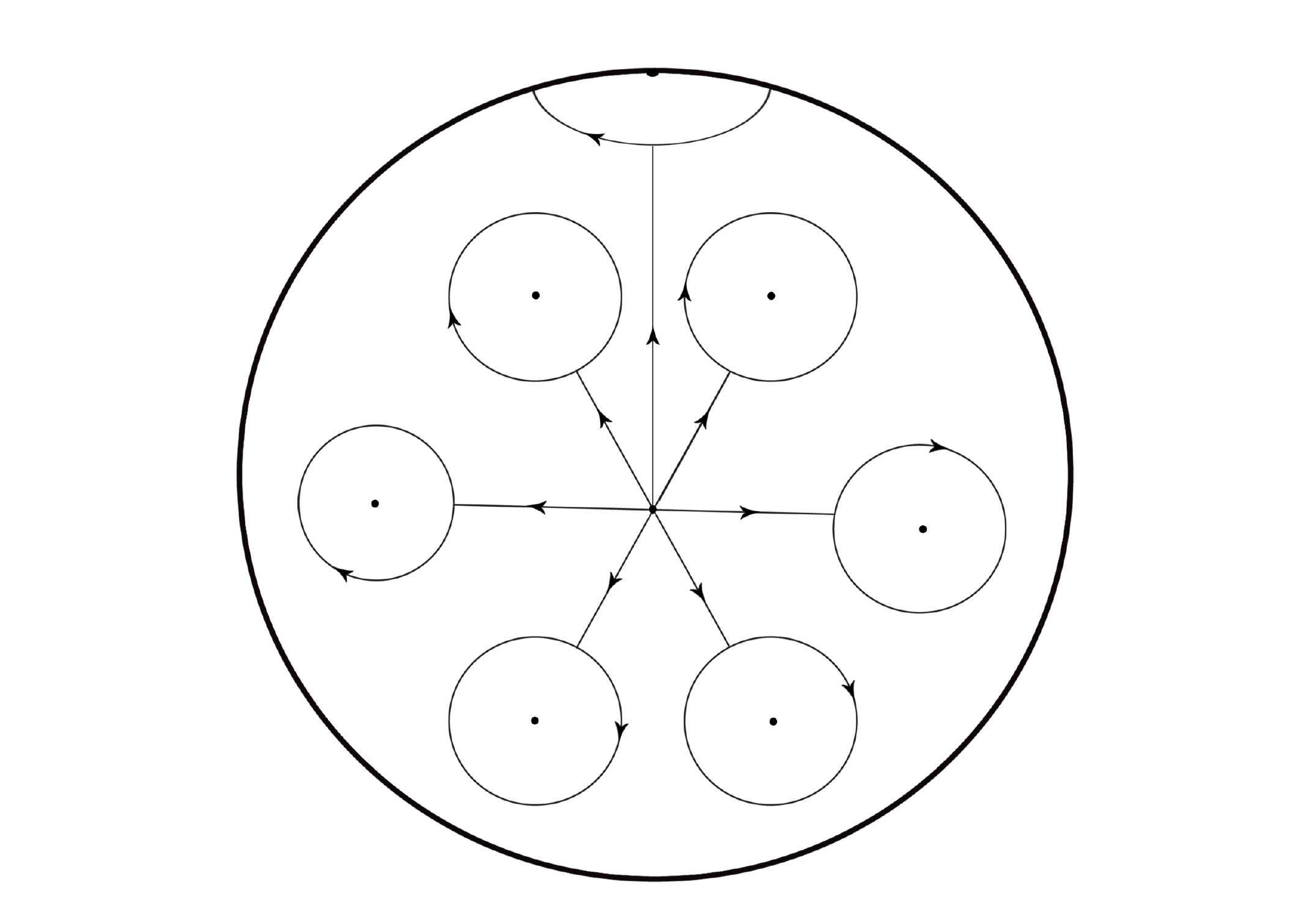}}}
\setbox0=\hbox{$\nwarrow$} 
\put(-173,139){$\nwarrow$\lower \ht0\hbox{$\delta_1$} }
\put(-161, 107){$\delta_2\uparrow$ }  
\setbox0=\hbox{$\nearrow$} 
\setbox0=\hbox{$\searrow$}
  \put(-190, 205){$\leftarrow\delta_\infty$} 
\put(-133, 170){$\leftarrow\sigma_1^*$} 
\put(-162,95){\lower \ht0\hbox{$\sigma_2^*$}$\nearrow$}   
\put(-265, 153){\lower \ht0\hbox{$\sigma_6^*$}$\nearrow$ } 
\put(-239, 217){\lower \ht0\hbox{$\sigma_\infty^*$}$\nearrow$ } 
\setbox0=\hbox{$\nearrow$} 
\setbox0=\hbox{$\leftarrow$} 
\put(-185,107){$0$}  
\put(-153, 174){$z_1$} 
\put(-110, 111){$z_2$}  
\put(-155, 56){$z_3$} 
\put(-220, 54){$z_4$} 
\put(-266, 115){$z_5$} 
\put(-220, 174){$z_6$} 
\put(-187, 233){$\infty$} 
$$
\end{figure}
\end{appendix} 

\providecommand{\bysame}{\leavevmode\hbox to3em{\hrulefill}\thinspace}
\providecommand{\MR}{\relax\ifhmode\unskip\space\fi MR }
\providecommand{\MRhref}[2]{%
   \href{http://www.ams.org/mathscinet-getitem?mr=#1}{#2}
}
\providecommand{\href}[2]{#2}


\end{document}